\documentclass[twocolumn]{autart}    

\usepackage{graphicx}
\usepackage{graphics}         	 	
\usepackage{comment}
\usepackage{color}
\usepackage{epstopdf}
\usepackage{amsmath,amsfonts,amssymb} 
\newcommand{\real}{\mathbb{R}}
\newcommand*\fvec[1]{\ensuremath{\mathbf{#1}}}                                	
\newcommand{\mc}{\mathcal}
\newcommand{\realnonnegative}{\mathbb{R}_{\geq0}}
\renewcommand{\circle}{\mathbb{S}^1}	
\newcommand{\torus}{\mathbb{T}}
\newcommand{\setdef}[2]{\{#1 \; | \; #2\}}
\newcommand{\map}[3]{#1: #2 \rightarrow #3}
\newcommand{\until}[1]{\{1,\dots, #1\}}

\newcommand{\subscr}[2]{{#1}_{\textup{#2}}}

\newcommand\oprocendsymbol{\hbox{$\square$}}
\newcommand\oprocend{\relax\ifmmode\else\unskip\hfill\fi\oprocendsymbol}

\usepackage{booktabs}
\usepackage{mathrsfs}

\renewcommand{\thefootnote}{\fnsymbol{footnote}}
\graphicspath{{./Figures/}}
\newcommand{\myvspace}[1]{}

\begin{document}

\begin{frontmatter}

\title{Synchronization and Power Sharing \\ for Droop-Controlled Inverters in Islanded Microgrids\thanksref{footnoteinfo}} 

\thanks[footnoteinfo]{Corresponding author J. W. Simpson-Porco.}

\author{John W. Simpson-Porco}\ead{johnwsimpsonporco@engineering.ucsb.edu},    
\author{Florian D\"{o}rfler}\ead{dorfler@engineering.ucsb.edu},               
\author[UCSB]{Francesco Bullo}\ead{bullo@engineering.ucsb.edu}  

\address[UCSB]{Center for Control, Dynamical Systems, and Computation, University of California, Santa Barbara, Santa Barbara, CA 93106, USA}  


\begin{abstract}                          
Motivated by the recent and growing interest in smart grid technology, we study the operation of DC/AC inverters in an inductive microgrid.
We show that a network of loads and DC/AC inverters equipped with power-frequency droop controllers can be cast as a Kuramoto model of phase-coupled oscillators.
This novel description, together with results from the theory of coupled oscillators, allows us to characterize the behavior of the network of inverters and loads.
Specifically, we provide a necessary and sufficient condition for the existence of a synchronized solution that is unique and locally exponentially stable.
We present a selection of controller gains leading to a desirable sharing of power among the inverters, and specify the set of loads which can be serviced without violating given actuation constraints.
Moreover, we propose a distributed integral controller based on averaging algorithms, which dynamically regulates the system frequency in the presence of a time-varying load. Remarkably, this distributed-averaging integral controller has the additional property that it preserves the power sharing properties of the primary droop controller.
Our results hold without assumptions on identical line characteristics or voltage magnitudes.

\end{abstract}

\end{frontmatter}
\section{Introduction}

A \emph{microgrid} is a low-voltage electrical network, heterogeneously
composed of distributed generation, storage, load, and managed autonomously
from the larger primary network. Microgrids are able to connect to the wide
area electric power system (WAEPS) through a Point of Common Coupling
(PCC), but are also able to ``island'' themselves and operate
independently. Energy generation within a microgrid can be highly
heterogeneous, any many these sources generate either variable frequency AC
power (wind) or DC power (solar), and are therefore interfaced with a
synchronous AC microgrid via power electronic devices called DC/AC (or
AC/AC) \emph{power converters}, or simply \emph{inverters}. In islanded
operation, inverters are operated as \emph{voltage sourced inverters}
(VSIs), which act much like ideal voltage sources. It is through these VSIs
that actions must be taken to ensure synchronization, security, power
balance and load sharing in the network.

\begin{figure}[t]
\begin{center}
\includegraphics[width=6.5cm]{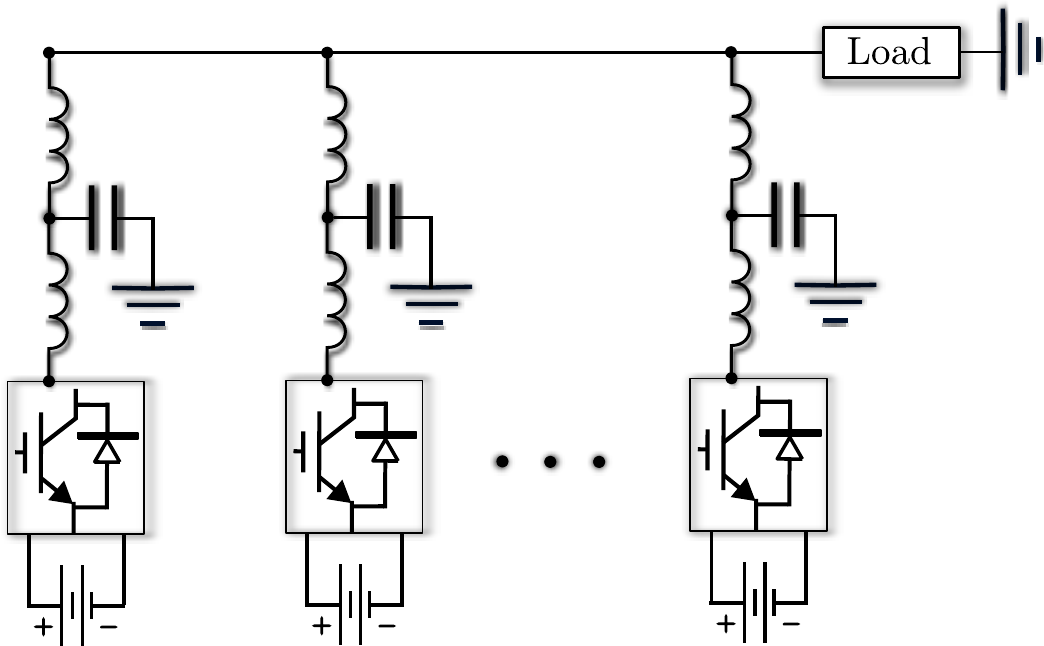}
\caption{Schematic of inverters operating in parallel.
}
\label{Fig:InvNet}
\end{center}  
\end{figure}

\textbf{\emph{Literature Review: }} A key topic of interest within the microgrid community is that of accurately sharing both active and reactive power among a bank of inverters operated \emph{in parallel}. Such a network is depicted in Figure \ref{Fig:InvNet}, in which each inverter transmits power directly to a common load.
Although several control architectures have been proposed to solve this problem, the so-called ``droop'' controllers have attracted the most attention, as they are ostensibly decentralized. The original reference for this methodology is \cite{MCC-DMD-RA:93}, where Chandorkar \emph{et. al.} introduce what we will refer to as the \emph{conventional droop controller}. 
For inductive lines, the droop controller attempts to emulate the behavior of a classical synchronous generator by imposing an inverse relation at each inverter between frequency and active power injection
\cite{PK:94}. Under other network conditions, the controller takes different forms \cite{JMG-LG-JM-MC-JM:05,YW-MC-JM-JMG-QMZ:11,QCZ-TH:13}.
Some representative references for the basic methodology are \cite{AT-HJ-TU-KM:97,SB-MC-PP-DP:02,JAPL-CLM-AGM:06,RM-AG-GL-FZ:08,YWL-CNK:09}
and \cite{JMG-JCV-JM-MC-LGDV:09}. Small-signal stability analyses for two inverters operating in parallel are presented under various assumptions in
\cite{EAAC-PCC-PFDG:02,MD-MNM-JWJ-AK:04,MNM-JJW-AK:07,YM-EFE-S:08}
and the references therein.
The recent work \cite{QCZ:13} highlights some drawbacks of the conventional droop method. Distributed controllers based on tools from synchronous generator theory and multi-agent systems have also been proposed for synchronization and power sharing. See \cite{ROS-JAF-RMM:07,WR-RWB-EMA:07} for a broad overview, and \cite{QCZ-TH:13,LABT-JPH-JM:11,SB-SZ:11,HX-ZQ-JS-AM:11} for various works.
%
%
%
%

Another set of literature relevant to our investigation is that pertaining to synchronization of phase-coupled oscillators, in particular the classic and celebrated \emph{Kuramoto model}. A generalization of this model considers $n \geq 2$ coupled oscillators, each represented by a phase $\theta_i \in \circle$ (the unit circle) and a natural frequency $\Omega_i \in \real$. The system of coupled oscillators obeys the dynamics
\myvspace{-0.5em}
\begin{equation}\label{Eq:Kuramoto}
\myvspace{-0.5em}
D_i\dot{\theta}_i = \Omega_i - \sum_{j=1}^n \nolimits a_{ij}\sin(\theta_i-\theta_j)\,,\,\, i \in \until{n},
\end{equation}
where $a_{ij} \geq 0$ is the coupling strength between the oscillators $i$ and $j$ and $D_i$ is the time constant of the $i^{th}$ oscillator.
\begin{figure}[t]
\begin{center}
\includegraphics[width=4cm]{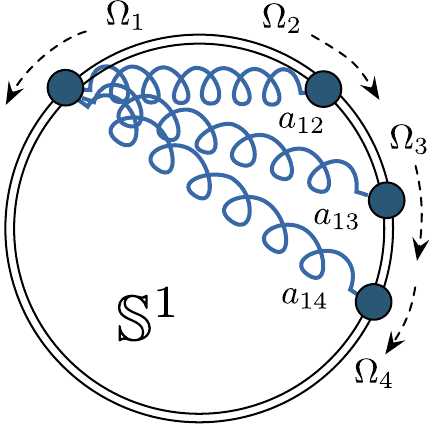}
\caption{Mechanical analog of a Kuramoto oscillator network. The particles have no inertia and do not collide with another.}
\label{Fig:MechAna}
\end{center}  
\end{figure}
%
Figure \ref{Fig:MechAna} shows a mechanical analog of \eqref{Eq:Kuramoto}, in which the oscillators can be visualized as a group of $n$ kinematic particles, constrained to rotate around the unit circle. The
particles rotate with preferred directions and speeds specified
by the natural frequencies $\Omega_i$, and are connected together by
elastic springs of stiffness $a_{ij}$. The rich dynamic behavior of
the system \eqref{Eq:Kuramoto} arises from the competition between the
tendency of each oscillator to align with its natural frequency $\Omega_i$,
and the synchronization enforcing coupling $a_{ij}\sin(\theta_i-\theta_j)$
with its neighbors. 
We refer to the recent
surveys
\cite{AA-ADG-JK-YM-CZ:08,SHS:00,FD-FB:10w}
for applications and theoretic results.

\myvspace{-0.8em}
\textbf{\emph{The Frequency-Droop Method: }} The \emph{frequency-droop method} constitutes one half of the conventional droop method. For inductive lines, the controller balances the active power demands in the network by instantaneously changing the frequency $\omega_i$ of the voltage signal at the $i^{th}$ inverter according to
\begin{equation}\label{Eq:Droop}
\omega_i = \omega^* - n_i(\subscr{P}{e,$i$}-P_i^*),
\end{equation}
where $\omega^*$ is a rated frequency, $\subscr{P}{e,$i$}$ is the active electrical power injection at bus $i$, and $P_i^*$ is the nominal active power injection. The parameter $n_i > 0$ is referred to as the \emph{droop coefficient}. 

\textbf{\emph{Limitations of the Literature:}} Despite forming the foundation for the operation of parallel VSIs, the frequency-droop control law \eqref{Eq:Droop} has never been subject to a nonlinear analysis \cite{QCZ:13}. No conditions have been presented under which the controller \eqref{Eq:Droop} leads the network to a synchronous steady state, nor have any statements been made about the convergence rate to such a steady state should one exist. Stability results that are presented rely on linearization for the special case of two inverters, and sometimes come packaged with extraneous assumptions \cite{RM-AG-GL-FZ:08,JMG-JCV-JM-MC-LGDV:09}. No guarantees are given in terms of performance. 
Schemes for power sharing based on ideas from multi-agent systems often deal directly with coordinating the real and reactive power injections of the distributed generators, and assume implicitly that a low level controller is bridging the gap between the true network physics and the desired power injections. 
Moreover, conventional schemes for frequency restoration typically rely on a combination of local integral action and separation of time scales, and are generally unable to maintain an appropriate sharing of power among the inverters (see Sections \ref{Sec:Select}--\ref{Sec:Secondary}).

\textbf{\emph{{Contributions: }}}The contributions of this paper are four-fold. 
First, we begin with our key observation that the equations governing a microgrid under the frequency-droop controller can be equivalently cast as a generalized Kuramoto model of the form \eqref{Eq:Kuramoto}. We present a necessary and sufficient condition for the existence of a locally exponentially stable and unique synchronized solution of the closed-loop, and provide a lower bound on the exponential convergence rate. We also state a robustified version of our stability condition which relaxes the assumption of fixed voltage magnitudes and admittances.
Second, we show rigorously\textemdash{}and without assumptions on large output impedances or identical voltage magnitudes\textemdash{}that if the droop coefficients are selected proportionally, then power is shared among the units proportionally. We provide explicit bounds on the set of 
serviceable loads.
Third, we propose a distributed ``secondary'' integral controller for frequency stabilization. Through the use of a distributed-averaging algorithm, the proposed controller dynamically regulates the network frequency to a nominal value, while preserving the proportional power sharing properties of the frequency-droop controller. We show that this controller is locally stabilizing, without relying on the classic assumption of a time-scale separation between the droop and integral control loops.
Fourth and finally, all results presented extend past the classic case of a parallel topology of inverters and hold for generic acyclic interconnections of inverters and~loads.

\textbf{\emph{Paper Organization: }}
The remainder of this section introduces some notation and reviews some fundamental material from algebraic graph theory, power systems and coupled oscillator theory. In Section \ref{Sec:ProbSetup} we motivate the mathematical models used throughout the rest of the work. In Section \ref{Sec:NonlinDroop} we perform a nonlinear stability analysis of the frequency-droop controller. Section \ref{Sec:Select} details results on power sharing and steady state bounds on power injections. In Section \ref{Sec:Secondary} we present and analyze our distributed-averaging integral controller. Finally, Section \ref{Section: Conclusions} concludes the paper and presents directions for future work.

\textbf{\emph{Preliminaries and Notation: }}\newline
\emph{Sets, vectors and functions}: Given a finite set $\mathcal{V}$, let $|\mathcal{V}|$ denote its cardinality. 
Given an index set $\mathcal{I}$ and a real valued 1D-array
$\{x_1,\ldots,x_{|\mathcal{I}|}\}$,
$\mathrm{diag}(\{x_i\}_{i\in\mathcal{I}}) \in
\real^{|\mathcal{I}|\times|\mathcal{I}|}$ is the associated diagonal
matrix. We denote the $n\times n$ identity matrix by $I_n$. Let
$\boldsymbol{1}_n$ and $\boldsymbol{0}_n$ be the $n$-dimensional vectors of
all ones and all zeros. For $z \in \real^n$, define $z^\perp \triangleq
\setdef{x \in \real^n}{z^Tx = 0}$ and $\boldsymbol{\sin}(z) \triangleq
(\sin(z_1),\ldots,\sin(z_n))^T \in \real^n$.

\emph{Algebraic graph theory: } 
We denote by $G(\mathcal{V},\mathcal{E},A)$ an undirected and weighted graph, where $\mathcal{V}$ is the set of nodes, $\mathcal{E} \subseteq \mathcal{V} \times \mathcal{V}$ is the set of edges, and $A \in \real^{|\mathcal{V}| \times |\mathcal{V}|}$ is the \emph{adjacency matrix}.
If a number $\ell \in \until{|\mathcal{E}|}$ and an arbitrary direction is assigned to each edge $\{i,j\} \in \mathcal{E}$, the \emph{node-edge incidence matrix} $B \in \real^{|\mathcal{V}|\times|\mathcal{E}|}$ is defined component-wise as $B_{k\ell} = 1$ if node $k$ is the sink node of edge $\ell$ and as $B_{k\ell} = -1$ if node $k$ is the source node of edge $\ell$, with all other elements being zero. 
For $x \in \real^{|\mathcal{V}|}$, $B^Tx \in \real^{|\mathcal{E}|}$ is the vector with components $x_i-x_j$, with $\{i,j\} \in \mathcal{E}$. 
%
If $\mathrm{diag}(\{a_{ij}\}_{\{i,j\}\in\mathcal{E}}) \in \real^{|\mathcal{E}| \times |\mathcal{E}|}$ is the diagonal matrix of edge weights, then the \emph{Laplacian matrix} is given by $L = B\mathrm{diag}(\{a_{ij}\}_{\{i,j\}\in\mathcal{E}})B^T$. 
If the graph is connected, then $\mathrm{ker}(B^T)=\mathrm{ker}(L) =
\mathrm{span}(\boldsymbol{1}_{|\mathcal{V}|})$, and $\mathrm{ker}(B) =
\emptyset$\@ for acyclic graphs. In this case, for every $x \in
\boldsymbol{1}_{|\mathcal{V}|}^\perp$, that is, $\sum_{i \in
  \mathcal{V}}x_i = 0$, there exists a unique $\xi \in
\real^{|\mathcal{E}|}$ satisfying Kirchoff's Current Law (KCL) $x = B\xi$
\cite{NB:97,LOC-CAD-EDK:87}. The vector $x$ is interpreted as nodal
injections, with $\xi$ being the associated edge flows. The Laplacian
matrix $L$ is positive semidefinite with eigenvalues $0 = \lambda_1(L) <
\lambda_2(L) \leq \cdots \leq \lambda_{|\mathcal{V}|}(L)$. We denote the
\emph{Moore-Penrose inverse} of $L$ by $L^\dagger$, and we recall
from~\cite{FD-FB:11d} the identity $LL^\dagger = L^\dagger L =
I_{|\mathcal{V}|} -
\frac{1}{|\mathcal{V}|}\boldsymbol{1}_{|\mathcal{V}|}\boldsymbol{1}_{|\mathcal{V}|}^T$.


\emph{Geometry on the $n$-torus}: The set $\circle$ denotes the \emph{unit circle}, an \emph{angle} is a point $\theta \in \circle$, and an \emph{arc} is a connected subset of $\circle$. With a slight abuse of notation, let $|\theta_1 - \theta_2|$ denote the \emph{geodesic distance} between two angles $\theta_1,\theta_2 \in \circle$. The \emph{$n$-torus} $\torus^n = \circle \times \cdots \times \circle$ is the Cartesian product of $n$ unit circles. For $\gamma \in [0,\pi/2[$ and a given graph $G(\mathcal{V},\mathcal{E},\cdot)$, let $\Delta_{G}(\gamma) =\{ \theta \in \mathbb T^{|\mathcal{V}|}:\, \max_{\{i,j\} \in \mathcal E} |\theta_{i} - \theta_{j}| \leq \gamma \}$ be the closed set of angle arrays $\theta = ( \theta_{1},\dots,\theta_{n})$ with neighboring angles $\theta_{i}$ and $\theta_{j}$, $\{i,j\} \in \mathcal E$ no further than $\gamma$ apart. 

\myvspace{-0.8em}
\emph{Synchronization:} Consider the first order phase-coupled oscillator model \eqref{Eq:Kuramoto} defined on a graph $G (\mathcal{V},\mathcal{E},\cdot)$. A solution $\map{\theta}{\realnonnegative}{\torus^{|\mathcal{V}|}}$ of \eqref{Eq:Kuramoto} is said to be \emph{synchronized} if (a) there exists a constant $\omega_{\mathrm{sync}} \in \real$ such that for each $t \geq 0$, $\dot{\theta}(t) = \omega_{\mathrm{sync}}\boldsymbol{1}_{|\mathcal{V}|}$ and (b) there exists a $\gamma \in [0,\pi/2[$ such that $\theta(t) \in \Delta_G(\gamma)$ for each $t \geq 0$. 

\myvspace{-0.8em}
\emph{AC Power Flow:} Consider a synchronous AC electrical network with $n$ nodes, purely inductive admittance matrix $Y \in j\mathbb{R}^{n\times n}$, nodal voltage magnitudes $E_i > 0$, and nodal voltage phase angles $\theta_i \in \circle$. The active electrical power $\subscr{P}{e,$i$} \in \real$ injected into the network at node $i \in \until{n}$ is given by \cite{PK:94}
\myvspace{-0.5em}
\begin{equation}\label{Eq:PowerFlow}
\subscr{P}{e,$i$} = \sum_{j=1}^{n}\nolimits E_iE_j|Y_{ij}|\sin(\theta_i-\theta_j).
\end{equation}

\myvspace{-2em}
\section{Problem Setup for Microgrid Analysis}
\myvspace{-0.8em}
\label{Sec:ProbSetup}\emph{Inverter Modeling: }
The standard approximation in the microgrid literature\textemdash{}and the one we adopt hereafter\textemdash{}
is to model an inverter as a \emph{controlled voltage source behind a reactance}.
This model is widely adopted among experimentalists in the microgrid field. Further modeling explanation can be found in \cite{ECF-LAA-LABT:08,QCZ-TH:13,BWW:92} and the references therein.
%

\begin{figure}[t]
\begin{center}
\includegraphics[width=6cm]{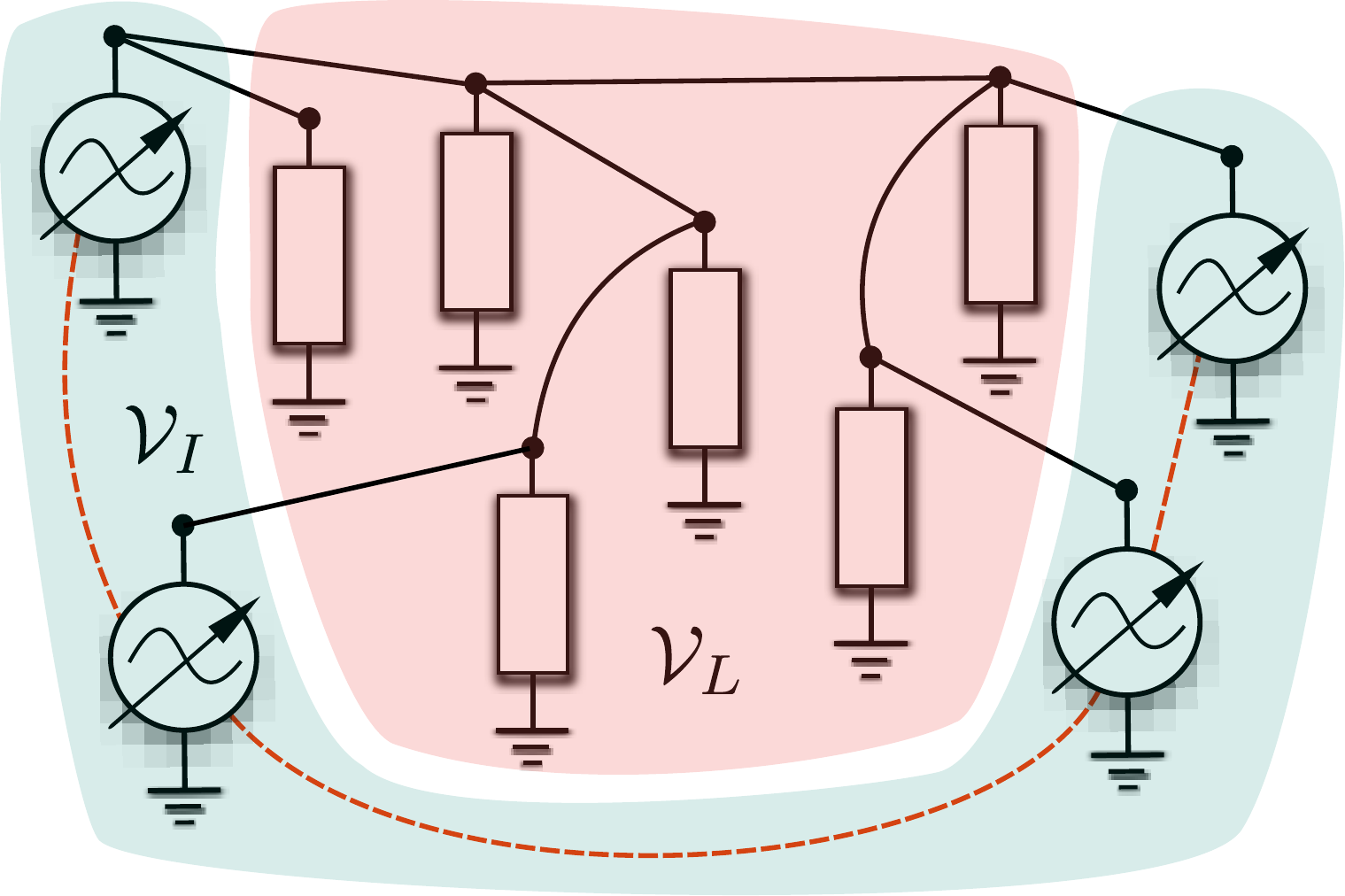}
\caption{Schematic illustration of a microgrid, with four inverters (nodes $\mathcal{V}_I$) supplying six loads (nodes $\mathcal{V}_L$) through an acyclic interconnection. The dotted lines between inverters represent communication links, which will be used exclusively in Section \ref{Sec:Secondary}.}
\label{Fig:InvTreeNet}
\end{center}
\end{figure}

\myvspace{-0.8em}
\emph{Islanded Microgrid Modeling: } Figure \ref{Fig:InvTreeNet} depicts an islanded microgrid containing both inverters and loads. Such an interconnection could arise by design, or spontaneously in a distribution network after an islanding event. An appropriate model is that of a weighted graph $G (\mathcal{V},\mathcal{E},A)$ with $|\mathcal{V}| = n$ nodes. 
%
We consider the case of inductive lines, and denote by $Y \in j\real^{n\times n}$ the bus admittance matrix of the 
network.\footnote{In some applications the inverter output impedance can be controlled to be highly inductive and dominate over any resistive effects in the network \cite{JMG-LG-JM-MC-JM:05}. In others, the control law \eqref{Eq:Droop} is inappropriate; see \cite{JMG-LG-JM-MC-JM:05,YW-MC-JM-JMG-QMZ:11,QCZ-TH:13} and Section \ref{Section: Conclusions}.}
We partition the set of nodes as $\mathcal{V} = \{\mathcal{V}_L,\mathcal{V}_I\}$, corresponding to loads and inverters.
For $\{i,j\} \in \mathcal{E}$, $Y_{ij}$ is the admittance of the edge between nodes $i$ and $j$. The output impedance of the inverter can be controlled to be purely susceptive, and we absorb its value into the line susceptances $-\mathrm{Im}(Y_{ij}) < 0$, $\{i,j\} \in \mathcal{E}$. 
To each node $i \in \until{n}$ we assign a harmonic voltage signal of the form $E_i(t) = E_i\cos(\omega^*t+\theta_i)$, where $\omega^* > 0$ is the nominal angular frequency, $E_i > 0$ is the voltage amplitude, and $\theta_i \in \circle$ is the voltage phase angle. We assume each inverter has precise measurements of its rolling time-averaged active power injection $\subscr{P}{e,$i$}(t)$ and of its frequency $\omega_{i}(t)$, see \cite{ECF-LAA-LABT:08} for details regarding this estimation.
%
The active power injection of each inverter into the network is restricted to the interval\,$[0,\overline{P}_i]$ where $\overline{P}_i$ is the \emph{rating} of inverter $i \in \mathcal{V}_I$. For the special case of a \emph{parallel} interconnection of inverters, as in Figure \ref{Fig:InvNet}, we will let $\mathcal{V}_L = \{0\}$ and $\mathcal{V}_I = \{1,\ldots,n-1\}$.

\myvspace{-1em}
\section{Analysis of Frequency-Droop Control}\label{Sec:NonlinDroop}
\myvspace{-0.8em}
We now connect the frequency-droop controller \eqref{Eq:Droop} to a network of first-order phase-coupled oscillators of the form \eqref{Eq:Kuramoto}. We restrict our attention to active power flows, and assume the voltage magnitudes $E_i$ are fixed at every bus. To begin, note that by defining $D_i \triangleq n_i^{-1}$ and by writing $\omega_i = \omega^* + \dot{\theta}_i$, we can equivalently write the frequency-droop controller \eqref{Eq:Droop} as
\myvspace{-0.5em}
\begin{equation}\label{Eq:BasicKuraDroop}
\myvspace{-0.5em}
D_i\dot{\theta}_i = {P}_i^* - P_{\mathrm{e},i}\,,\quad i \in \mathcal{V}_I,
\end{equation}
where $P_i^* \in [0,\overline{P}_i]$ is a selected nominal value.\footnote[2]{We make no assumptions regarding the selection of droop coefficients. See Section \ref{Sec:Select} for more on choice of coefficients.} 
Note that $\dot{\theta}_i$ is the \emph{deviation} of the frequency at inverter $i$ from the nominal frequency $\omega^*$.
Using the active load flow equations \eqref{Eq:PowerFlow}, the droop controller \eqref{Eq:BasicKuraDroop} becomes
\myvspace{-0.5em}
\begin{equation}\label{Eq:KuraDroop}
\myvspace{-0.5em}
D_i\dot{\theta}_i = P_i^* - \sum_{j=1}^n\nolimits E_iE_j|Y_{ij}|\sin(\theta_i-\theta_j),\quad i \in \mathcal{V}_I.
\end{equation}
For a constant power loads $\{P_i^*\}_{i \in \mathcal{V}_L}$ we must also satisfy the $|\mathcal{V}_L|$ power balance equations 
\myvspace{-0.5em}
\begin{equation}\label{Eq:PowerBal}
\myvspace{-0.5em}
0 = P_i^* - \sum_{j=1}^n \nolimits E_iE_j|Y_{ij}|\sin(\theta_i-\theta_j)\,,\quad i \in \mathcal{V}_L\,.
\end{equation}
%
If due to failure or energy shortage an inverter $i \in \mathcal{V}_I$ is unable to supply power support to the network, we formally set $D_i = P_i^* = 0$, which reduces \eqref{Eq:KuraDroop} to a load as in \eqref{Eq:PowerBal}.
If the droop-controlled system \eqref{Eq:KuraDroop}--\eqref{Eq:PowerBal} reaches synchronization, then we can \textemdash{} without loss of generality \textemdash{} transform our coordinates to a \emph{rotating frame of reference}, where the synchronization frequency is zero and the study of synchronization reduces to the study of equilibria.
In this case, it is known that the equilibrium point of interest for the differential-algebraic system shares the same stability properties as the same equilibrium of the corresponding singularly perturbed system \cite[Theorem 13.1]{H-DC:11}, where the constant power loads $P_{i}^{*}$ are replaced by a frequency-dependent loads $P_{i}^{*} - D_{i} \dot \theta_{i}$, for $i \in \mathcal{V}_L$ and for some sufficiently small $D_{i} > 0$.
%
We can now identify a singularly perturbed droop-controlled system of the form \eqref{Eq:KuraDroop}--\eqref{Eq:PowerBal} with a network of Kuramoto oscillators described by \eqref{Eq:Kuramoto} and arrive at the following insightful relation.
\myvspace{-0.8em}
%
\begin{lem}\label{Lem:EquivOfModels}\textbf{\emph{(Equivalence of Perturbed Droop-Controlled System and Kuramoto Model).}}
The following two models are equivalent:
\myvspace{-0.8em}
\begin{enumerate}
\item[(i)] The singularly perturbed droop-controlled network \eqref{Eq:KuraDroop}--\eqref{Eq:PowerBal}, with frequency-dependent loads $P_{i}^{*} - D_{i} \dot \theta_{i}$ with $D_{i}> 0$ instead of constant power loads $P_i^* \in \real$ ($i \in \mathcal{V}_L$), droop coefficients $n_i = 1/D_{i} > 0$, nominal power injections $P_i^* \in \real$ ($i \in \mathcal{V}_I$), nodal voltage phases $\theta_i \in \circle$, nodal voltages magnitudes $E_i > 0$, and bus admittance matrix $Y \in j\real^{n \times n}$.
\item[(ii)] The generalized Kuramoto model \eqref{Eq:Kuramoto}, with time constants $D_i > 0$, natural frequencies $\Omega_i \in \real$, phase angles $\theta_i \in \circle$ and coupling weights $a_{ij} > 0$.
\end{enumerate}
\myvspace{-0.8em}
Moreover, the parametric quantities of the two models are related via $P_i^* = \Omega_i$  and $E_iE_j|Y_{ij}| = a_{ij}$.
\end{lem}
\myvspace{-1em}

In light of Lemma \ref{Lem:EquivOfModels} and for notational simplicity, we define the matrix of time constants (inverse droop coefficients) $D \triangleq \mathrm{diag}(\boldsymbol{0}_{|\mathcal{V}_L|},\{D_i\}_{i \in \mathcal{V}_I})$, the vector of loads and nominal power injections $P^* \triangleq (P_1^*,\ldots,P_{n}^*)^T$, and for $\{i,j\} \in \mathcal{E}$ we write $a_{ij} \triangleq E_iE_j|Y_{ij}|$. 
The drop-controlled system \eqref{Eq:KuraDroop}--\eqref{Eq:PowerBal} then reads in vector notation~as
\myvspace{-0.5em}
\begin{equation}\label{Eq:CombinedKuraDroop}
\myvspace{-0.5em}
D\dot{\theta} = P^* - B\mathrm{diag}(\{a_{ij}\}_{\{i,j\}\in\mathcal{E}})\boldsymbol{\sin}(B^T\theta),
\end{equation}
where $\theta \triangleq (\theta_1,\ldots,\theta_n)^T$ and $B \in \real^{n \times |\mathcal{E}|}$ is the node-edge incidence matrix of the underlying graph $G(\mathcal{V},\mathcal{E},A)$. 
A natural question now arises: under what conditions on the power injections, network topology, admittances, and droop coefficients does the differential-algebraic closed-loop system \eqref{Eq:KuraDroop}--\eqref{Eq:PowerBal} possess a stable, synchronous solution?
\smallskip

\myvspace{-0.8em}
\begin{thm}\label{Thm:Stab}\textbf{\emph{(Existence and Stability of Sync'd Solution).}} Consider the frequency-droop controlled system \eqref{Eq:KuraDroop}--\eqref{Eq:PowerBal} defined on an acyclic network
with node-edge incidence matrix $B$. Define the scaled power imbalance $\omega_{\rm avg}$ by $\omega_{\textup{avg}} \triangleq (\sum_{i=1}^nP_i^*) / (\sum_{i \in \mathcal{V}_I}D_i) \in \real$, and let $\xi \in \real^{|\mathcal{E}|}$ be the unique vector of edge power flows satisfying KCL, given implicitly by $P^*-\omega_{\rm avg}D\boldsymbol{1}_{n}=B\xi$\@. The following two statements are equivalent:
%
\myvspace{-0.5em}
\begin{enumerate}
\item[(i)] \emph{\textbf{Synchronization: }} There exists an arc length $\gamma \in [0,\pi/2[$ such that the closed-loop system \eqref{Eq:KuraDroop}--\eqref{Eq:PowerBal} possess a locally exponentially stable and unique synchronized solution $t \mapsto \theta^*(t) \in \Delta_G(\gamma)$ for all $t \geq 0$;
\item[(ii)] \emph{\textbf{Flow Feasibility: }}The power flow is feasible, i.e., 
\myvspace{-0.5em}
\begin{equation} \label{Eq:DorflerCondition}
\Gamma \triangleq \|\mathrm{diag}(\{a_{ij}\}_{\{i,j\}\in\mathcal{E}})^{-1}\xi\|_{\infty} < 1.
\end{equation}
\end{enumerate}
\myvspace{-0.8em}
If the equivalent statements (i) and (ii) hold true, then the quantities $\Gamma \in [0,1[$ and $\gamma \in [0,\pi/2[$ are related uniquely via $\Gamma = \sin(\gamma)$, and the following statements hold:
\myvspace{-0.8em}
\begin{enumerate}
\item[a)] \emph{\bf Explicit Synchronized Solution:} The synchronized solution satisfies $\theta^*(t) =  \theta_0 + \left(\omega_{\textup{sync}}t\boldsymbol{1}_{n} \right) \pmod{2 \pi}$ for some $\theta_0 \in \Delta_G(\gamma)$, where $\omega_{\textup{sync}} = \omega_{\rm avg}$, and the synchronized angular differences satisfy $\boldsymbol{\sin}(B^T\theta^*) = \mathrm{diag}(\{a_{ij}\}_{\{i,j\}\in\mathcal{E}})\xi$;
\item[b)] \emph{\textbf{Explicit Synchronization Rate: }} The local exponential synchronization rate is no worse than
\myvspace{-0.5em}
\begin{equation}
\myvspace{-0.5em}
\lambda \triangleq \frac{\lambda_2(L)}{\max_{i \in \mathcal{V}_I} D_i}\sqrt{1-\Gamma^2},
\end{equation}
where $L = B\mathrm{diag}(\{a_{ij}\}_{\{i,j\}\in\mathcal{E}})B^T$ is the Laplacian matrix of the network with weights~$\{a_{ij}\}_{\{i,j\}\in\mathcal{E}}$.
\end{enumerate}
\end{thm}
\myvspace{-0.5em}
\begin{rem}\label{Rem:PhysInterp}\emph{\textbf{(Physical Interpretation)}}
From the droop controller \eqref{Eq:BasicKuraDroop}, it holds that $P^* - \omega_{\rm sync}D\boldsymbol{1}_{n} \in \boldsymbol{1}_n^\perp$ is the vector of steady state power injections. The power injections therefore satisfy the Kirchoff current law, and $\xi \in \real^{|\mathcal{E}|}$ is the associated vector of power flows along edges \cite{LOC-CAD-EDK:87}. Physically, the parametric condition \eqref{Eq:DorflerCondition} therefore states that the active power flow along each edge be feasible, i.e., less than the physical maximum $a_{ij} = E_iE_j|Y_{ij}|$. While the necessity of this condition seems plausible, its sufficiency is perhaps surprising. Theorem \ref{Thm:Stab} shows that equilibrium power flows are invariant under constant scaling of all droop coefficients, as overall scaling of $D$ appears inversely in $\omega_{\mathrm{avg}}$. Although grid stress varies with specific application and loading, the condition \eqref{Eq:DorflerCondition} is typically satisfied with a large margin of safety -- a practical upper bound for $\gamma$ would be $10^{\circ}$.
%
\end{rem}
\myvspace{-1.5em}
\begin{pf}
%
%
To begin, note that if a solution $t \mapsto \theta(t)$ to the system
\eqref{Eq:CombinedKuraDroop} is frequency synchronized, then by definition there exists an
$\omega_{\rm sync} \in \real$ such that $\dot{\theta}(t) = \omega_{\rm
  sync}\boldsymbol{1}_{n}$ for all $t \geq 0$. Summing over all equations
\eqref{Eq:KuraDroop}--\eqref{Eq:PowerBal} gives $\omega_{\rm sync} =
\omega_{\rm avg}$. Without loss of generality, we can consider the
\emph{auxiliary system} associated with \eqref{Eq:CombinedKuraDroop}
defined~by
\myvspace{-0.5em}
\begin{equation}\label{Eq:Aux}
\myvspace{-0.5em}
D\dot{\theta} = \widetilde{P} - B\mathrm{diag}(\{a_{ij}\}_{\{i,j\}\in\mathcal{E}})\boldsymbol{\sin}(B^T\theta),
\end{equation}
where $\widetilde{P}_i = P_i^*$ for $i \in \mathcal{V}_L$ and $\widetilde{P}_i = P_i^* - \omega_\mathrm{avg}D_i$ for $i \in \mathcal{V}_I$. Since $\widetilde{P} \in \boldsymbol{1}_n^\perp$, system \eqref{Eq:Aux} has the property that $\widetilde{\omega}_{\mathrm{avg}} = 0$ and represents the dynamics \eqref{Eq:CombinedKuraDroop} in a reference frame rotating at an angular frequency $\omega_{\mathrm{avg}}$. Thus, frequency synchronized solutions of \eqref{Eq:CombinedKuraDroop} correspond one-to-one with \emph{equilibrium} points of the system \eqref{Eq:Aux}.
Given the Laplacian matrix $L = B\mathrm{diag}(\{a_{ij}\}_{\{i,j\}\in\mathcal{E}})B^T$, \eqref{Eq:Aux} can be equivalently rewritten in the insightful form
\myvspace{-0.5em}
\begin{equation}\label{Eq:Aux-2}
\myvspace{-0.5em}
D\dot{\theta} = B\mathrm{diag}(\{a_{ij}\}_{\{i,j\}\in\mathcal{E}}) \cdot \bigl( B^TL^{\dagger}\widetilde{P} - \boldsymbol{\sin}(B^T\theta) \bigr) \,.
\end{equation}
Here, we have made use of the facts that $LL^{\dagger} = L^{\dagger}L = I_n
- \frac{1}{n} \fvec 1_{n} \fvec 1_{n}^{T}$ and $\widetilde P \in
\boldsymbol{1}_n^\perp$.
Since $\mathrm{ker}(B) = \emptyset$, equilibria of \eqref{Eq:Aux-2} must
satisfy $B^TL^\dagger\widetilde{P} = B^TL^\dagger B \xi =
\boldsymbol{\sin}(B^T\theta)$. We claim that $B^TL^{\dagger}B =
\mathrm{diag}(\{a_{ij}\}_{\{i,j\}\in\mathcal{E}})^{-1}$. To see this, define $X \triangleq B^TL^\dagger B$ and notice that $X\mathrm{diag}(\{a_{ij}\}_{\{i,j\}\in\mathcal{E}})B^T = B^TL^\dagger(B\mathrm{diag}(\{a_{ij}\}_{\{i,j\}\in\mathcal{E}})B^T) = B^TL^\dagger L = B^T$. Since $\mathrm{ker}(B) = \emptyset$, it therefore holds that $X\mathrm{diag}(\{a_{ij}\}_{\{i,j\}\in\mathcal{E}}) = I_{|\mathcal{E}|}$ and the result follows. Hence, equilibria of \eqref{Eq:Aux-2} satisfy
\begin{equation}\label{Eq:ReducedDorfler}
\mathrm{diag}(\{a_{ij}\}_{\{i,j\}\in\mathcal{E}})^{-1}\xi = \boldsymbol{\sin}(B^T\theta).
\end{equation}
Equation \eqref{Eq:ReducedDorfler} is uniquely solvable for $\theta^* \in
\Delta(\gamma)$, $\gamma \in {[0,\pi/2[}$, if and only if
$\Gamma \triangleq \max_{\{i,j\}\in\mathcal{E}} (\xi_{ij}/a_{ij}) \leq
\sin(\gamma)$.
Since the right-hand side of the condition $\Gamma \leq \sin(\gamma)$ is a concave and monotonically increasing function of $\gamma \in [0,\pi/2[$, there exists an equilibrium $\theta^{*} \in \Delta_{G}(\gamma)$ for some $\gamma \in [0,\pi/2[$ if and only if the condition $\Gamma \leq \sin(\gamma)$ is true with the strict inequality sign for $\gamma = \pi/2$. This leads immediately to the claimed condition $\Gamma < 1$. 
In this case, the explicit equilibrium angles are then obtained from the $n$ decoupled equations \eqref{Eq:ReducedDorfler}. See \cite[Theorems 1, 2(G1)]{fd-mc-fb:11v-pnas} for additional information.
Local exponential stability of the equilibrium $\theta^* \in \Delta(\gamma)$ is established by recalling the equivalence between the index-1 differential-algebraic system \eqref{Eq:Aux-2} and an associated reduced set of pure differential equations (see also the proof of (b)). In summary, the above discussion shows the equivalence of (i) and (ii) and statement~(a).
%
%

To show statement (b), consider the linearization of the dynamics \eqref{Eq:Aux} about the equilibrium $\theta^* \in \Delta(\gamma)$ given by
\begin{equation*}
\frac{\mathrm{d}}{\mathrm{d}t}\begin{bmatrix}\boldsymbol{0}_{|\mathcal{V}_L|}\\ \Delta\theta_I\end{bmatrix} = -\begin{bmatrix}I_{|\mathcal{V}_L|} & \boldsymbol{0} \\ \boldsymbol{0} & D_{I}^{-1}\end{bmatrix}\begin{bmatrix}L_{LL} & L_{LI} \\ L_{IL} & L_{II}\end{bmatrix}\begin{bmatrix}\Delta \theta_L\\ \Delta\theta_I\end{bmatrix},
\end{equation*}
where we have partitioned the matrix $L(\theta^*) = B\mathrm{diag}(\{a_{ij}\cos(\theta_i^*-\theta_j^*)\}_{\{i,j\}\in\mathcal{E}})B^T$ according to load nodes $\mathcal{V}_L$ and inverter nodes $\mathcal{V}_I$, and defined $D_I \triangleq \mathrm{diag}(\{D_i\}_{i \in \mathcal{V}_I})$. Since $\theta^* \in \Delta_G(\gamma)$, the matrix $L(\theta^*)$ is a Laplacian and thus positive semidefinite with a simple eigenvalue at zero corresponding to rotational invariance of the dynamics under a uniform shift of all angles. 
It can be easily verified that the upper left block $L_{LL}$ of $L(\theta^{*})$ is nonsingular  \cite{FD-FB:11d}, or equivalently, $\theta^{*}$ is a {\it regular} equilibrium point.
Solving the set of $|\mathcal{V}_L|$ algebraic equations and substituting into the dynamics for $\Delta\theta_I$, we obtain $\mathrm{d}(\Delta\theta_I)/\mathrm{d}t = -D_{I}^{-1}L_{\rm red}(\theta^*)\Delta \theta_I$, where $L_{\rm red} \triangleq L_{II} - L_{IL}L_{LL}^{-1}L_{LI}$. The matrix $L_{\rm red}(\theta^*) \in \real^{|\mathcal{V}_I|\times|\mathcal{V}_I|}$ is also a Laplacian matrix, and therefore shares the same properties as $L(\theta^*)$ \cite{FD-FB:11d}.
%
Thus, it is the \emph{second} smallest eigenvalue of $D_{I}^{-1}L_{\rm
  red}(\theta^*)$ which bounds the convergence rate of the linearization,
and hence the local convergence rate of the dynamics \eqref{Eq:Aux}.
A simple bound on $\lambda_2(D_{I}^{-1}L_{\rm red}(\theta^*))$ can be obtained via the Courant-Fischer Theorem \cite{CDM:01}.
For $x \in \boldsymbol{1}_{|\mathcal{V}_I|}^\perp$, let $y = D_I^{1/2}x$, and note that
 $x^TL_{\rm red}(\theta^*)x/(x^TD_Ix)=y^T\!D_I^{-1/2}\!L_{\rm red}(\theta^*)D_I^{-1/2}y/(y^T\!y)$.%
\,Thus,\,$y \!\in\! (D_{I}^{-1/2}\boldsymbol{1}_{|\mathcal{V}_I|})^\perp$ is an eigenvector of $D_I^{-1/2}L_{\rm red}(\theta^*)D_I^{-1/2}$ with eigenvalue~$\mu \in \real$ if and only if $x = D_{I}^{-1/2}y$ is an eigenvector of $D_{I}^{-1}L_{\rm red}(\theta^*)$ with eigenvalue $\mu$. For $y \neq \boldsymbol{0}_{|\mathcal{V}_I|}$, we obtain
\begin{align*}
&\lambda_{2}(D_I^{-1}L_{\rm red}(\theta^*)) = \!\!\min_{y \in (D_{I}^{-1/2}\boldsymbol{1}_{|\mathcal{V}_I|})^\perp}\!\! \frac{y^TD_I^{-\frac{1}{2}}L_{\rm red}(\theta^*)D_I^{-\frac{1}{2}}y}{y^Ty}\\ 
&= \!\min_{x \in \boldsymbol{1}_{|\mathcal{V}_I|}^\perp}\!\! \frac{x^TL_{\rm red}(\theta^*)x}{x^TD_Ix}
 \geq \frac{1}{\max_{i \in \mathcal{V}_I}D_i} \!\min_{x \in \boldsymbol{1}_{|\mathcal{V}_I|}^\perp}\!\! \frac{x^TL_{\rm red}(\theta^*)x}{x^Tx}\\
&\geq \frac{\lambda_2(L_{\rm red}(\theta^*))}{\max_{i \in \mathcal{V}_I}D_i} \geq \frac{\lambda_2(L(\theta^*))}{\max_{i \in \mathcal{V}_I}D_i},
\end{align*}
where we have made use of the \emph{spectral interlacing} property of Schur complements \cite{FD-FB:11d} in the final inequality. Since $\theta^* \in \Delta_G(\gamma)$, the eigenvalue $\lambda_2(L(\theta^*))$ can be further bounded as $\lambda_2(L(\theta^*)) \geq \lambda_2(L)\cos(\gamma)$, where $L = B\mathrm{diag}(\{a_{ij}\}_{\{i,j\}\in\mathcal{E}})B^T$ is the Laplacian with~weights $\{a_{ij}\}_{\{i,j\}\in\mathcal{E}}$. This fact and the identity $\cos(\gamma) = \cos(\sin^{-1}(\Gamma)) = \sqrt{1-\Gamma^2}$ complete the proof. 
\oprocend \end{pf}
\myvspace{-2em}
An analogous stability result for inverters operating in parallel now follows as a corollary. 
\myvspace{-0.8em}
\begin{cor}\label{Thm:StabParallel}\emph{\textbf{(Existence and Stability of Sync'd Solution for Parallel Inverters).}}
Consider a parallel interconnection of inverters, as depicted in Figure \ref{Fig:InvNet}. The following two statements are equivalent:
\myvspace{-0.8em}
\begin{enumerate}
\item[(i)] \emph{\textbf{Synchronization: }} There exists an arc length $\gamma \in [0,\pi/2[$ such that the closed-loop system \eqref{Eq:CombinedKuraDroop} possess a locally exponentially stable and unique synchronized solution $t \mapsto \theta^*(t) \in \Delta_G(\gamma)$ for all $t \geq 0$;
\item[(ii)] \emph{\textbf{Power Injection Feasibility: }} 
\begin{equation} \label{Eq:DorflerConditionParallel}
\Gamma \triangleq \max_{i\in\mathcal{V}_I} |(P_i^*-\omega_{\mathrm{avg}}D_{i})/a_{i0}| < 1.
\end{equation}
\end{enumerate}
\end{cor}
\myvspace{-1em}
\begin{pf}
For the parallel topology of Figure \ref{Fig:InvNet} there is one load fed by $n-1$ inverters, and the incidence matrix of the graph $G(\mathcal{V},\mathcal{E},A)$ takes the form $B = \left[ -\boldsymbol{1}_{n-1} \,\, I_{n-1}
\right]^T$. Letting $\widetilde{P}$ be as in the previous proof, we note that $\xi$ is given uniquely as $\xi =  (B^TB)^{-1}B^T\widetilde{P}$. 
In this case, a set of straightforward but tedious matrix calculations reduce condition \eqref{Eq:DorflerCondition} to condition \eqref{Eq:DorflerConditionParallel}.
\oprocend
\end{pf}

\myvspace{-1em}
Our analysis so far has been based on the assumption that each term $a_{ij} \triangleq E_iE_j|Y_{ij}|$ is a {\em constant and known} parameter for all $\{i,j\} \in \mc E$. In a realistic power system, both effective line susceptances magnitudes and voltage magnitudes are dynamically adjusted by additional controllers.
%
%
Our analysis so far has been based on the assumption that each term $a_{ij} \triangleq E_iE_j|Y_{ij}|$ is a {\em constant and known} parameter for all $\{i,j\} \in \mc E$. In a realistic power system, both effective line susceptances magnitudes and voltage magnitudes are dynamically adjusted by additional controllers. The following result states that as long as these controllers can regulate the effective susceptances and nodal voltages above prespecified lower bounds $|\underline{Y_{ij}}|$ and $\underline{E_j}$, the stability results of Theorem \ref{Thm:Stab} go through with little modification.
\begin{cor}\emph{\textbf{(Robustified Stability Condition).}}
\label{Cor:Robust}
Consider the frequency-droop controlled system \eqref{Eq:KuraDroop}--\eqref{Eq:PowerBal}. 
Assume that the nodal voltage magnitudes satisfy $E_i > \underline{E_i} >0$ for all $i \in \until{n}$, and that the line susceptance magnitudes satisfy $|Y_{ij}| \geq |\underline{Y_{ij}}| > 0$ for all $\{i,j\} \in \mc E$. For $\{i,j\} \in \mathcal{E}$, define $\underline{a_{ij}} \triangleq \underline{E_{i}}\underline{E_j} \underline{|Y_{ij}|}$. 
The following two statements are equivalent:
\myvspace{-0.5em}
\begin{enumerate}
\item[(i)] \emph{\textbf{Robust Synchronization: }} For all possible voltage magnitudes $E_i > \underline{E_i}$ and line susceptance magnitudes $|Y_{ij}| \geq |\underline{Y_{ij}}|$, there exists an arc length $\gamma \in [0,\pi/2[$ such that the closed-loop system \eqref{Eq:KuraDroop}--\eqref{Eq:PowerBal} possess a locally exponentially stable and unique synchronized solution $t \mapsto \theta^*(t) \in \Delta_G(\gamma)$ for all $t \geq 0$; and
%
\item[(ii)] \emph{\textbf{Worst Case Flow Feasibility: }} The active power flow is feasible for the worst case voltage and line susceptances magnitudes, that is,
\myvspace{-0.5em}
\begin{equation}\label{Eq:RobustDorfler}
\myvspace{-0.5em}
\|\mathrm{diag}(\{\underline{a_{ij}}\}_{\{i,j\}\in\mathcal{E}})^{-1}\xi\|_{\infty} < 1.
\end{equation}
\end{enumerate}
\end{cor}
\myvspace{-0.9em}
\begin{pf}
The result follows by noting that $a_{ij}$ (resp. $\underline{a}_{ij}$) appears exclusively in the \emph{denominator} of \eqref{Eq:DorflerCondition} (resp. \eqref{Eq:RobustDorfler}), and that the vector 
$\xi \in \real^{|\mathcal{E}|}$ defined in Theorem \ref{Thm:Stab} does not depend on the voltages or line susceptances.\oprocend \end{pf}
Finally, regarding the assumption of purely inductive lines, we note that since the eigenvalues of a matrix are continuous functions of its entries, the exponential stability property established in Theorem \ref{Thm:Stab} is robust, and the stable synchronous solution persists in the presence of sufficiently small line conductances \cite{HDC-CCC:95}. This robustness towards lossy lines is also illustrated in the simulation study of Section \ref{Sec:Sim}.
\myvspace{-1.5em}
\section{Power Sharing and Actuation Constraints}\label{Sec:Select}
\myvspace{-0.8em}
While Theorem \ref{Thm:Stab} gives the necessary and sufficient condition for the existence of a synchronized solution to the closed-loop system \eqref{Eq:KuraDroop}--\eqref{Eq:PowerBal}, it offers no immediate guidance on how to select the control parameters $P_i^*$ and $D_i$ to satisfy the actuation constraint $\subscr{P}{e,$i$} \in [0,\overline{P}_i]$.
The following definition gives the proper criteria for selection. 


\myvspace{-0.8em}
\begin{defn}\textbf{\emph{(Proportional Droop Coefficients).}}\label{Def:Propor}
The droop coefficients are selected proportionally if $P_i^*/D_i =
P_j^*/D_j$ and $P_i^*/\overline{P}_i = P_j^*/\overline{P}_j$ for all $i,j
\in \mathcal{V}_I$.
\end{defn}

\myvspace{-0.8em}
\begin{thm}\label{Thm:PowerFlowConstraintsTree}\textbf{\emph{(Power Flow Constraints and Power Sharing).}}
Consider a synchronized solution of the frequency-droop controlled system \eqref{Eq:KuraDroop}--\eqref{Eq:PowerBal}, and let the droop coefficients be selected proportionally. Define the total load $P_L \triangleq \sum_{i \in \mathcal{V}_L}P_i^*$. The following two statements are equivalent:
\myvspace{-0.8em}
\begin{enumerate}
\item[(i)] \emph{\textbf{Injection Constraints: }} $0 \leq \subscr{P}{e,$i$} \leq \overline{P}_i$, $\,\,\,\forall i \in \mathcal{V}_I$;
\item[(ii)] \emph{\textbf{Load Constraint: }} $-\sum_{j \in \mathcal{V}_I}\overline{P}_j \leq P_L \leq 0.$
\end{enumerate}
\myvspace{-0.8em}
Moreover, the inverters share the total load $P_L$ proportionally according to their power ratings, that is, $\subscr{P}{e,$i$}/\overline{P}_i = P_{\mathrm{e},j}/\overline{P}_j$, for each $i,j \in \mathcal{V}_I$.
\end{thm}
\myvspace{-2em}
\begin{pf}
From \eqref{Eq:BasicKuraDroop}, the steady state active power injection at
each inverter is given by $\subscr{P}{e,$i$} = P_i^* - \omega_{\rm
  sync}D_i$. By imposing $\subscr{P}{e,$i$} \geq 0$ for each $i \in
\mathcal{V}_I$, substituting the expression for $\omega_{\rm sync}$, and
rearranging terms, we obtain, for each $i \in \mathcal{V}_I$,
\myvspace{-0.8em}
\begin{align*}
\myvspace{-0.8em}
\subscr{P}{e,$i$} = P_i^* &- \left(\frac{P_L + \sum_{j\in\mathcal{V}_I}P_j^*}{\sum_{j\in\mathcal{V}_I}D_j}\right)D_i \geq 0 \\ \Longleftrightarrow\,\,\, &P_L \leq -\sum_{j\in\mathcal{V}_I} \left(P_j^* - \frac{P_i^*}{D_i}D_j\right) = 0,
\end{align*}
where in the final equality we have used Definition \ref{Def:Propor}. Along with the observation that $P_{\mathrm{e},i} \geq 0$ if and only if $P_{e,j} \geq 0$ ($i,j \in \mathcal{V}_I$), this suffices to show that $0 \leq \subscr{P}{e,$i$}$ for each $i \in \mathcal{V}_I$ if and only if $P_L \leq 0$. If we now impose for $i \in \mathcal{V}_I$ that $\subscr{P}{e,$i$} \leq \overline{P}_i$ and again use the expression for $\omega_{\rm sync}$ along with Definition \ref{Def:Propor}, a similar calculation~yields
\myvspace{-0.8em}
\begin{align*}
\myvspace{-0.8em}
\subscr{P}{e,$i$} \leq \overline{P}_i &\Longleftrightarrow P_L \geq
-\frac{\overline{P}_i}{P_i^*}\sum_{j\in\mathcal{V}_I} \nolimits P_j^* =- \sum_{j \in \mathcal{V}_I}\nolimits\overline{P}_j.
\end{align*}
Along with the observation that $P_{\mathrm{e},i} \leq \overline{P}_i$ if and only if $P_{\mathrm{e},j} \leq \overline{P}_j$ ($i,j \in \mathcal{V}_I$), this shows that $\subscr{P}{e,$i$} \leq \overline{P}_i$ for each $i \in \mathcal{V}_I$ if and only if the total load $P_L$ satisfies the above inequality. In summary, we have demonstrated two if and only if inequalities, which when taken together show the equivalence of $(i)$ and $(ii)$. To show the final statement, note that the fraction of the rated power capacity injected by the $i^{th}$ inverter is given by
\myvspace{-0.8em}
$$
\myvspace{-0.8em}
\frac{\subscr{P}{e,$i$}}{\overline{P}_i} = \frac{P_i^*-\omega_{\rm sync}D_i}{\overline{P}_i}
=  \frac{P_j^*-\omega_{\rm sync}D_j}{\overline{P}_j} = \frac{P_{\mathrm{e},j}}{\overline{P}_j},
$$
for each $j \in \mathcal{V}_I$. This completes the proof. \oprocend
\end{pf}
\myvspace{-1.5em}
Power sharing results for parallel inverters supplying a single load follow as a corollary of Theorem \ref{Thm:PowerFlowConstraintsTree}, with $P_L = P_0^*$. 
Note that the coefficients $D_i$ must be selected with global knowledge. The droop method therefore requires a centralized design for power sharing despite its decentralized implementation. We remark that Theorem \ref{Thm:PowerFlowConstraintsTree} holds independently of the network voltage magnitudes and line admittances. 
%
%
\myvspace{-1.5em}
\section{Distributed PI Control in Microgrids}\label{Sec:Secondary}
\myvspace{-0.8em}
As is evident from the expression for $\omega_{\rm sync}$ in Theorem \ref{Thm:Stab}, the frequency-droop method typically leads to a deviation of the steady state operating frequency $\omega^* + \omega_{\rm sync}$ from the nominal value $\omega^*$. Again in light of Theorem \ref{Thm:Stab}, it is clear that modifying the nominal active power injection $P_i^*$ via the transformation $P_i^* \longrightarrow P_i^* - \omega_{\mathrm{sync}}D_i$ (for $i\in\mathcal{V}_I$) in the controller \eqref{Eq:KuraDroop} will yield zero steady state frequency deviation (c.f. the auxiliary system \eqref{Eq:Aux} with $\widetilde{\omega}_{\rm sync} = 0$). Unfortunately, the information to calculate $\omega_{\rm sync}$ is not available locally at each inverter.
As originally proposed in \cite{MCC-DMD-RA:93}, after the frequency of each inverter has converged to $\omega_{\mathrm{sync}}$, a slower, ``secondary'' control loop can be used locally at each inverter. Each local secondary controller slowly modifies the nominal power injection $P_i^*$ until the network frequency deviation is zero. This procedure implicitly assumes that the measured frequency value $\dot{\theta}_i(t)$ is a good approximation of $\omega_{\mathrm{sync}}$, and relies on a separation of time-scales between the fast, synchronization-enforcing primary droop controller and the slower secondary integral controller. This methodology is employed in \cite{MCC-DMD-RA:93,RL-PP:00,JMG-JCV-JM-MC-LGDV:09}. For large droop coefficients $D_i$, this approach can be particularly slow (Theorem \ref{Thm:Stab} (b)), with this slow response leading to an inability of the method to dynamically regulate the network frequency in the presence of a time-varying load. Moreover, these decentralized integral controllers destroy the power sharing properties established by the primary droop controller.
\begin{figure}[t]
\begin{center}
\includegraphics[width=8cm]{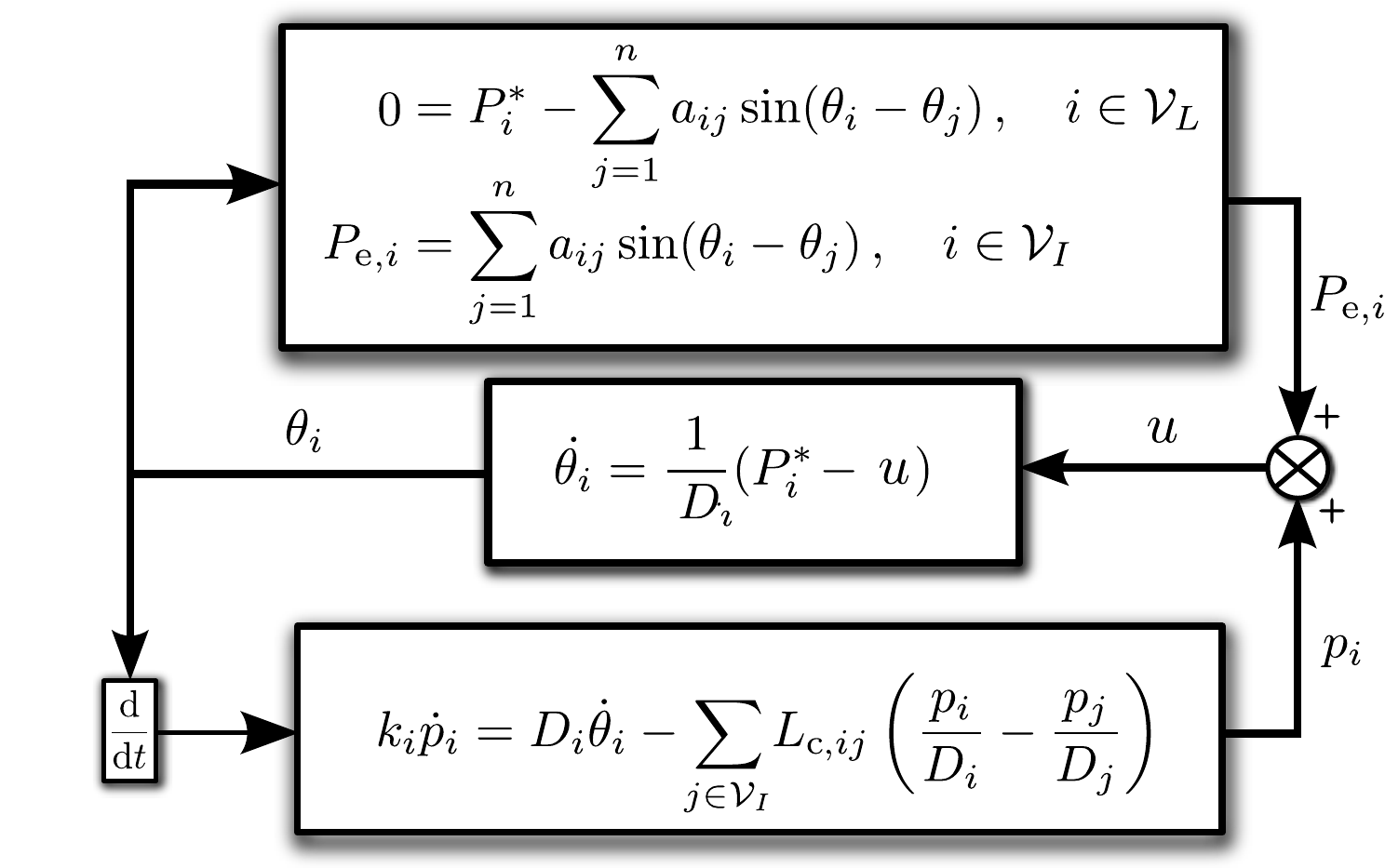}
\caption{Feedback diagram for the DAPI controller. By itself, the upper control loop is the droop controller \eqref{Eq:BasicKuraDroop}. 
}
\label{Fig:DroopFeedback3}
\end{center}  
\end{figure}

%
%
In what follows, we pursue an alternative scheme for frequency restoration
which does not implicitly rely on a separation of time-scales as in
\cite{MCC-DMD-RA:93,RL-PP:00,JMG-JCV-JM-MC-LGDV:09}. Assuming the existence of a communication network among the inverters, we expand on the
conventional frequency-droop design \eqref{Eq:Droop} and propose the
\emph{distributed-averaging proportional-integral (DAPI) controller}
\myvspace{-0.8em}
\begin{align}
\myvspace{-0.8em}
	D_i\dot{\theta}_i &= {P}_i^* - p_{i} - \subscr{P}{e,$i$} 
	\,,\quad i \in \mathcal{V}_I\,,
	\label{eq: primary control}
	\\
	 k_i \dot p_{i} &=D_{i}\dot \theta_{i}-\sum_{j \in \mathcal{V}_I} L_{\mathrm{c},ij}\left(\frac{p_i}{D_i}-\frac{p_j}{D_j}\right),\,
	\, i \in \mathcal{V}_I \,,
	\label{eq: secondary control}
\end{align}
where $p_i \in \real$ is an auxiliary power variable
and $k_{i} > 0$ is a gain, for each $i \in \mathcal{V}_I$.\footnote[2]{The presented results extend to discrete time and asynchronous communication, see \cite{FB-JC-SM:09}.} The matrix $L_{\rm c} \in \real^{|\mathcal{V}_I|\times|\mathcal{V}_I|}$ is the Laplacian matrix corresponding to a weighted, undirected and connected \emph{communication graph} $G_{\rm c} (\mathcal{V}_{I},\mathcal{E}_{\rm c},A_{\rm c})$ between the inverters, see Figure \ref{Fig:InvTreeNet}. 
The DAPI controller \eqref{eq: primary control}--\eqref{eq: secondary control} is depicted in Figure~\ref{Fig:DroopFeedback3}, and will be shown to have the following two key properties. First of all, the controller is able to quickly regulate the network frequency under large and rapid variations in load. Secondly, the controller accomplishes this regulation while \emph{preserving} the power sharing properties of the primary droop controller \eqref{Eq:BasicKuraDroop}. 
The closed-loop dynamics resulting from the DAPI controller \eqref{eq: primary control}--\eqref{eq: secondary control} are given by
\myvspace{-0.8em}
\begin{align}
	0 =&\; P_{i}^* - \sum_{j=1}^n \nolimits a_{ij} \sin(\theta_{i} - \theta_{j})\,, \quad i \in \mathcal{V}_L\,,
	\label{eq: load -- closed loop}
		\\
	D_{i} \dot \theta_{i} =&\;  {P}_i^* - p_{i} -  \sum_{j=1}^n \nolimits a_{ij} \sin(\theta_{i} - \theta_{j})
	\,,\,\,  i \in \mathcal{V}_I\,,
	\label{eq: primary control -- closed loop}
	\\
	k_i\dot p_{i} =&\; {P}_i^* - p_{i} -  \sum_{j=1}^n \nolimits a_{ij} \sin(\theta_{i} - \theta_{j})\nonumber\\ &\;-\sum_{j\in\mathcal{V}_I} \nolimits L_{\mathrm{ c},ij}\left(\frac{p_i}{D_i}-\frac{p_j}{D_j}\right) 
	\,, \quad i \in \mathcal{V}_I \,.
	\label{eq: secondary control -- closed loop}
\end{align}
The following theorem (see Appendix \ref{App:AdaptiveProof} for the proof) establishes local stability of the desired equilibrium of \eqref{eq: load -- closed loop}--\eqref{eq: secondary control -- closed loop} as well as the power sharing properties of the DAPI controller.

\myvspace{-0.8em}
\begin{thm}\emph{\textbf{(Stability of DAPI-Controlled Network).}}
\label{Theorem: Stability of adaptive droop controller}
Consider an acyclic network of droop-controlled inverters and loads in which the inverters can communicate through the weighted graph $G_c$, as described by the closed-loop system \eqref{eq: load -- closed loop}--\eqref{eq: secondary control -- closed loop} with parameters $P_{i}^{*} \in [0,\overline{P}_i], D_{i}>0$ and $k_{i}>0$ for $i \in \mathcal{V}_I$, and connected communication Laplacian $L_{\rm c} \in \real^{|\mathcal{V}_I|\times |\mathcal{V}_I|}$.
The following two statements are equivalent:
\myvspace{-0.8em}
\begin{enumerate}
\item[(i)] \emph{\textbf{Stability of Droop Controller: }} The droop control stability condition \eqref{Eq:DorflerCondition} holds;
\item[(ii)] \emph{\textbf{Stability of DAPI Controller: }} There exists an arc length $\gamma \in [0,\pi/2[$ such that the system \eqref{eq: primary control -- closed loop}--\eqref{eq: secondary control -- closed loop} possess a locally exponentially stable and unique equilibrium $\bigl( \theta^{*} , p^{*} \bigr) \in \Delta_G(\gamma) \times \real^{|\mathcal{V}_I|}$.
\end{enumerate}
\myvspace{-0.8em}
If the equivalent statements (i) and (ii) hold true, then the unique equilibrium is given as in Theorem \ref{Thm:Stab} $(ii)$, along with $p_{i}^{*} = D_{i} \subscr{\omega}{avg}$ for $i \in \mathcal{V}_I$. Moreover, if the droop coefficients are selected proportionally, then the DAPI controller preserves the proportional power sharing property of the primary droop controller.
\end{thm}
\myvspace{-0.8em}
%
\myvspace{-1.2em}
Note that Theorem \ref{Theorem: Stability of adaptive droop controller} asserts the exponential stability of an \emph{equilibrium} of the closed-loop \eqref{eq: load -- closed loop}--\eqref{eq: secondary control -- closed loop}, and hence, a synchronization frequency $\omega_{\rm sync}$ of zero. The network therefore synchronizes to the nominal frequency $\omega^*$\@.
\myvspace{-1em}
\section{Simulation Study}\label{Sec:Sim}

\begin{figure}[t]
\begin{center}
\includegraphics[width=1.0\columnwidth]{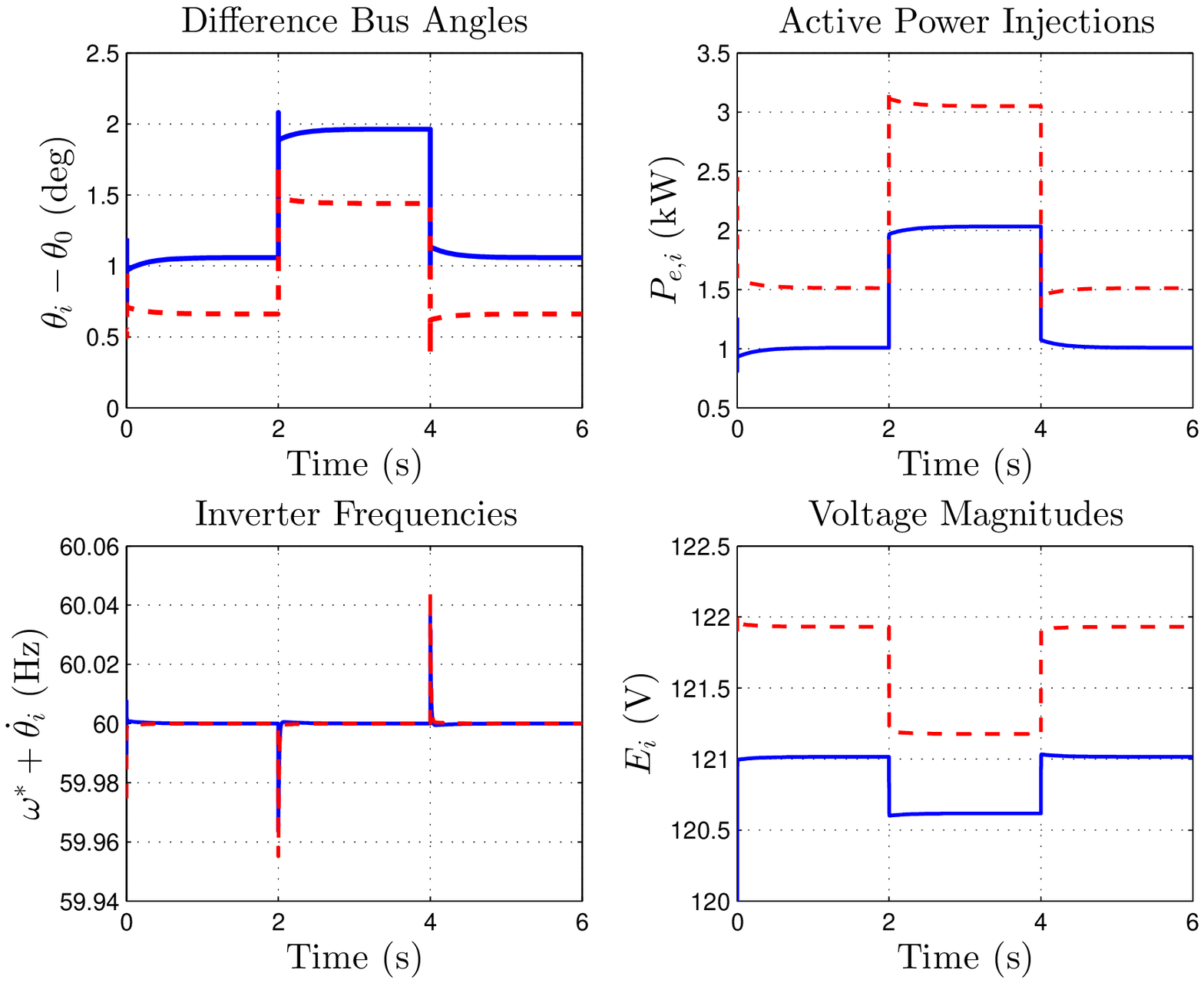}
\caption{DAPI controlled closed-loop \eqref{eq: load -- closed loop}--\eqref{eq: secondary control -- closed loop} for two inverters supplying a load which changes at $t \in {\{2\textup{s},4\textup{s}\}}$. }
\label{Fig:TwoInvSim}
\end{center}  
\end{figure}

\myvspace{-0.8em}
%
We now illustrate the performance of our proposed  DAPI controller \eqref{eq: primary control}--\eqref{eq: secondary control} and its robustness to unmodeled voltage dynamics (see Corollary \ref{Cor:Robust}) and lossy lines in a simulation scenario. We consider two inverters operating in parallel and supplying a variable load. The voltage magnitude at each inverter is controlled via the \emph{voltage-droop} method
\myvspace{-0.6em}
\begin{equation}\label{Eq:QVDroop}
\myvspace{-0.6em}
E_i = E_i^* - m_i\left(Q_{\mathrm{e},i}-Q_i^*\right)\,,\quad i \in \{1,2\},
\end{equation}
where $E_i^* > 0$ (resp. $Q_i^* \in \real$) is the nominal voltage magnitude (resp. nominal reactive power injection) at inverter $i \in \{1,2\}$, $m_i > 0$ is the voltage-droop coefficient, and $Q_{\mathrm{e},i} \in \real$ is the reactive power injection (see \cite{PK:94} for details on reactive power).
%
The simulation parameters are reported in Table \ref{Tab:TwoInvSimPar}.
Note the effectiveness of the proposed DAPI controller \eqref{eq: primary control}--\eqref{eq: secondary control} in quickly regulating the system frequency. The spikes in local frequency displayed in Figure  \ref{Fig:TwoInvSim} (c) are due to the rapid change in load, and can be further suppressed by increasing the gains $k_i$. This additional degree of freedom allows for a selection of primary droop coefficients $D_i$ much smaller than is typical in the literature ($10^3$W$\cdot$s, compared to roughly $10^5$W$\cdot$s), allowing the power injections (Figure \ref{Fig:TwoInvSim} (b)) to respond quickly during transients. 

\begin{table}
\begin{center}
\caption{Parameter values for simulation in Figure \ref{Fig:TwoInvSim}.} The choice of resistances corresponds to a resistance/reactance ratio of one half.
{\renewcommand{\arraystretch}{1}
\begin{tabular}{llr}
\toprule
Parameter & Symbol & Value \\
\midrule
Nom. Frequency & $\omega^*/2\pi$ & 60\,Hz\\
Nom. Voltages & $E_i^*$ & [120,\,122]\,V\\
Output/Line Induc. & $L_i$ & [0.7,\,0.5]\,mH\\
Output/Line Resist. & $R_i$ & [0.14,\,0.1]\,$\Omega$\\
Inv. Ratings ($P$) & $P_i^* = \overline{P}_i$ & [2,\,3]\,kW \\
Inv. Ratings ($Q$) & $Q_i^*$ & [1,\,1]\,kvar \\
Load ($P$) & $P_0^*(t)$ & $P_0^* \in \{-2.5,-5\}$kW\\
Load ($Q$) & $Q_0^*(t)$ & $Q_0^* \in \{-1,-2\}$kvar\\
$\omega$--Droop Coeff. & $D_i$ & [4,\,6] $\times 10^3$ $\mathrm{W\cdot s}$\\
$E$--Droop Coef. & $m_i$ & [1,\,1] $\times 10^{-3}$ $\frac{\rm V}{\rm var}$\\
Sec. Droop Coeff. & $k_i$ & $10^{-6}$\,s\\
Comm. Graph & $G_{\rm comm}$ & Two nodes, one edge\\
Comm. Laplacian & $L_{\rm c}$ & $(1000\, \mathrm{W}\mathrm{s})\cdot \begin{bmatrix}1 & -1 \\ -1 & 1\end{bmatrix}$ \\
\bottomrule
\end{tabular}
}
\label{Tab:TwoInvSimPar}
\end{center}
\end{table}

\myvspace{-0.8em}
\section{Conclusions}\label{Section: Conclusions}
\myvspace{-0.8em}
We have examined the problems of synchronization, power sharing, and secondary control among droop-controlled inverters by leveraging tools from the theory of coupled oscillators, along with ideas from classical power systems and multi-agent systems. An issue not addressed in this work is a nonlinear analysis of reactive power sharing, as an analysis of the voltage-droop method \eqref{Eq:QVDroop} which yields simple and physically meaningful algebraic conditions for the existence of a solution is difficult to perform. Moreover, in the case of strongly mixed line conditions $\mathrm{Im}(Y_{ij}) \simeq \mathrm{Re}(Y_{ij})$, neither of the control laws \eqref{Eq:Droop} nor \eqref{Eq:QVDroop} are appropriate. A provably functional control strategy for general interconnections and line conditions is an open and exciting problem.
\myvspace{-1em}
\begin{ack}                               
\myvspace{-0.8em}
This work was supported in part by the National Science Foundation NSF CNS-1135819 and by the National Science and Engineering Research Council of Canada. We wish to thank H. Bouattour, Q.-C. Zhong and J. M. Guerrero for their insightful comments and suggestions.
\end{ack}
\myvspace{-0.8em}
\bibliographystyle{plain}
\bibliography{alias,Main,FB}
\appendix
\myvspace{-0.8em}
\section{Proof of Theorem \ref{Theorem: Stability of adaptive droop controller}}
\label{App:AdaptiveProof}
\myvspace{-0.8em}

Consider the closed-loop \eqref{eq: load -- closed loop}--\eqref{eq: secondary control -- closed loop} arising from the DAPI controller \eqref{eq: primary control}--\eqref{eq: secondary control}. We formulate our problem in the error coordinates $\widetilde p_{i}(t) \triangleq p_{i}(t) - D_{i}\subscr{\omega}{avg}$, 
and write \eqref{eq: load -- closed loop}--\eqref{eq: secondary control -- closed loop} in vector notation~as 
\begin{align}
	\mathcal{D}\begin{bmatrix}\boldsymbol{0}_{|\mathcal{V}_L|}\\ \dot{\theta}_I\end{bmatrix} &= \widetilde{P} - P_{\rm e} - \begin{bmatrix}\boldsymbol{0}_{|\mathcal{V}_L|}\\\widetilde{p}\end{bmatrix},\label{eq: load -- closed loop -- rot2}
	\\
       K\dot{\widetilde{p}} = \widetilde{P}_I &- P_{\mathrm{e},I}-(I_{|\mathcal{V}_I|}+L_{\rm c}D_I^{-1})\widetilde{p},\label{eq: secondary control -- closed loop -- rot2}
\end{align}
where we have defined $P_{\rm e} \triangleq B\mathcal{A}\boldsymbol{\sin}(B^T\theta) = (P_{\mathrm{e},L},P_{\mathrm{e},I})$, $\mathcal{A} \triangleq  \mathrm{diag}(\{a_{ij}\}_{\{i,j\}\in\mathcal{E}})$, $\mathcal{D} \triangleq \mathrm{blkdiag}(I_{|\mathcal{V}_L|},D_I)$, $K \triangleq \mathrm{diag}(\{k_i\}_{i \in \mathcal{V}_I})$ and partitioned the vector of power injections by load and inverter nodes as  $\widetilde{P} = (\widetilde{P}_L,\widetilde{P}_I)^T$.
Equilibria of \eqref{eq: load -- closed loop -- rot2}--\eqref{eq: secondary control -- closed loop -- rot2} satisfy
\begin{equation}\label{Eq:EquilibriaEqn}
\boldsymbol{0} = \underbrace{\begin{bmatrix}I_{|\mathcal{V}_L|} & \boldsymbol{0} & \boldsymbol{0} \\ \boldsymbol{0} & D_I^{-1} & I_{|\mathcal{V}_I|} \\ \boldsymbol{0} & I_{|\mathcal{V}_I|} & D_I + L_{\rm c}\end{bmatrix}}_{Q_1}\underbrace{\begin{bmatrix}I_{|\mathcal{V}_L|} & \boldsymbol{0} & \boldsymbol{0} \\ \boldsymbol{0} & I_{|\mathcal{V}_I|} & \boldsymbol{0} \\ \boldsymbol{0} & \boldsymbol{0} & D_I^{-1}\end{bmatrix}}_{Q_2}\underbrace{\begin{bmatrix}\widetilde{P} - P_{\mathrm{e}} \\ -\widetilde{p}\end{bmatrix}}_{x}.
\end{equation}
The positive semidefinite matrix $Q_1$ has one dimensional kernel spanned by $(\boldsymbol{0}_{|\mathcal{V}_L|},D_I\boldsymbol{1}_{|\mathcal{V}_I|},-\boldsymbol{1}_{|\mathcal{V}_I|})$, while $Q_2$ is positive definite. Note however that since $\widetilde{P} - P_{\rm e} \in \boldsymbol{1}_n^\perp$, and $Q_2x = (\widetilde{P}-P_{\rm e},-D_I^{-1}\widetilde{p})$, it holds that $Q_2x \notin \mathrm{ker}(Q_1)$. Thus, \eqref{Eq:EquilibriaEqn} holds if and only if $x = \boldsymbol{0}_{n+|\mathcal{V}_I|}$; that is, $\widetilde{p} = \widetilde{p}^* =  \boldsymbol{0}_{|\mathcal{V}_I|}$ and $\widetilde{P}-P_{\rm e} = \boldsymbol{0}_n$. Equivalently, from Theorem \ref{Thm:Stab}, the latter equation is solvable for a unique (modulo rotational symmetry) value $\theta^* \in \Delta_G(\gamma)$ if and only if the parametric condition \eqref{Eq:DorflerCondition} holds.

To establish the local exponential stability of the equilibrium $(\theta^*,\widetilde{p}^*)$, we linearize the DAE \eqref{eq: load -- closed loop -- rot2}--\eqref{eq: secondary control -- closed loop -- rot2} about the regular fixed point $(\theta^*,\widetilde{p}^*)$ and eliminate the resulting algebraic equations, as in the proof of Theorem~\ref{Thm:Stab}. The Jacobian $J(\theta^*,\widetilde{p}^*)$ of the reduced system of ordinary differential equations can then be factored as $J(\theta^*,\widetilde{p}^*) = -Z^{-1}X$, where $Z = \mathrm{blkdiag}(I_{|\mathcal{V}_I|},K)$ and 
$$
X = \underbrace{\begin{bmatrix} D_I^{-1} & I_{|\mathcal{V}_I|} \\ I_{|\mathcal{V}_I|} & L_{\rm c} + D_I\end{bmatrix}}_{=X_1 = X_1^T} \underbrace{\begin{bmatrix}L_{\rm red}(\theta^*) & \boldsymbol{0} \\ \boldsymbol{0} & D_I^{-1}\end{bmatrix}}_{=X_2 = X_2^T}.
$$
Thus, the problem of local exponential stability of $(\theta^*,\widetilde{p}^*)$ reduces to the generalized eigenvalue problem $-X_1X_2v = \lambda Z v$, where $\lambda \in \real$ is an eigenvalue and $v \in \real^{2|\mc V_{I}|}$ is the associated eigenvector.
We will proceed via a continuity-type argument. Consider momentarily a perturbed version of $X_1$, denoted by $X_1^{\epsilon}$, obtained by adding the matrix $\epsilon I_{|\mathcal{V}_I|}$ to the lower-right block of $X_1$, where $\epsilon \geq 0$. Then for every $\epsilon > 0$, $X_1^\epsilon$ is positive definite. Defining $y \triangleq Zv$, we can write the generalized eigenvalue problem $X_1^{\epsilon}X_2v = -\lambda Z v$ as $X_2Z^{-1}y = -\lambda (X_1^{\epsilon})^{-1}y$. The matrices on both left and right of this generalized eigenvalue problem are now symmetric, with $X_2Z^{-1} = \mathrm{blkdiag}(L_{\rm red},D_I^{-1}K^{-1})$ having a simple eigenvalue at zero corresponding to rotational symmetry. By applying the Courant-Fischer Theorem to this transformed problem, we conclude, for $\epsilon >0$ and modulo rotational symmetry, that all eigenvalues are real and negative.

Now, consider again the unperturbed case with $\epsilon = 0$. Notice that the matrix $X_2$ is positive semidefinite with kernel spanned by $(\boldsymbol{1}_{|\mathcal{V}_I|},\boldsymbol{0}_{|\mathcal{V}_I|})$ corresponding to rotational symmetry, while $X_1$ is positive semidefinite with kernel spanned by $(-D_I\boldsymbol{1}_{|\mathcal{V}_I|},\boldsymbol{1}_{|\mathcal{V}_I|})$. Since $\mathrm{image}(L_{\rm red}(\theta^*)) = \boldsymbol{1}^\perp_{|\mathcal{V}_I|}$, $\mathrm{image}(X_2) \cap \mathrm{ker}(X_1) = \{\boldsymbol{0}_{2|\mathcal{V}_I|}\}$,   that is, $X_2v$ is never in the kernel of $X_1$. Thus we conclude that $\mathrm{ker}(X_1X_2) = \mathrm{ker}(X_2)$. Now we return to the original eigenvalue problem in the form $X_1X_2v = -\lambda Z v$.
Since the eigenvalues of a matrix are continuous functions of the matrix entries, and $\mathrm{ker}(X_1X_2) = \mathrm{ker}(X_2)$, we conclude that the number of negative eigenvalues does not change as $\epsilon \rightarrow 0^+$, and the eigenvalues therefore remain real and negative. Hence, the equilibrium $(\theta^*,\widetilde{p}^*)$ of the  DAE \eqref{eq: load -- closed loop -- rot2}--\eqref{eq: secondary control -- closed loop -- rot2} is (again, modulo rotational symmetry) locally exponentially stable.

To show the final statement, note from the modified primary controller \eqref{eq: primary control} that the steady state power injection at inverter $i \in \mathcal{V}_I$ is given by $P_{\mathrm{e},i} = P_i^* - p_i(t = \infty) = P_i^* - \omega_{\rm avg}D_i$, which is exactly the steady state power injection when only the primary droop controller \eqref{Eq:BasicKuraDroop} is used. The result then follows from Theorem \eqref{Thm:PowerFlowConstraintsTree}. This completes
the proof of Theorem \ref{Theorem: Stability of adaptive droop controller}.
\oprocend

\end{document}